\documentclass{article}
\usepackage{spconf, amsmath, graphicx}

\usepackage{enumitem}
\setlist{nosep, leftmargin=14pt}

\usepackage{amsmath, amssymb, amsfonts, amsthm, mathtools, physics, mathrsfs, stmaryrd}
\usepackage{notations}
\usepackage{siunitx}
\usepackage{braket}
\usepackage{subcaption}

\usepackage{tikz}
\usetikzlibrary{spy}
\usepackage[linesnumbered,ruled,vlined]{algorithm2e}
\SetKwInput{KwData}{Input}
\SetKwInput{KwResult}{Ouput}

\title{Off-the-grid covariance-based super-resolution fluctuation microscopy}
\name{Bastien Laville $^{1}$, Laure Blanc-Féraud $^{1}$, Gilles Aubert $^{2}$.}
\address{%
    $^{1}$ Université Côte d’Azur, CNRS, Inria, I3S, Morpheme project, France.  \\
    $^{2}$Université Côte d’Azur, CNRS, LJAD, France. 
}

\usepackage{fancyhdr}

\fancypagestyle{FirstPage}{
\setlength{\headheight}{40pt}
\chead{\small \emph{© 2021 IEEE.  Personal use of this material is permitted. Permission from IEEE must be obtained for all other uses, in any current or future media, including reprinting/republishing this material for advertising or promotional purposes, creating new collective works, for resale or redistribution to servers or lists, or reuse of any copyrighted component of this work in other works.}}
}

\begin{document}
%

\thispagestyle{FirstPage}

\maketitle

\begin{abstract}
    Super-resolution fluorescence microscopy overcomes blurring arising from light diffraction, allowing the reconstruction of fine scale details in biological structures. Standard methods come at the expense of long acquisition time and/or harmful effects on the biological sample, which makes the problem quite challenging for the imaging of body cells.
    A promising new avenue is the exploitation of molecules fluctuations, allowing live-cell imaging with good spatio-temporal resolution through common microscopes and conventional fluorescent dyes.
    Several numerical algorithms have been developed in the literature and used for fluctuant time series. These techniques are developed within the discrete setting, namely the super-resolved image is defined on a finer grid than the observed images. On the contrary, gridless optimisation does not rely on a fine grid and is rather an optimisation of Dirac measures in number, amplitudes and positions. In this work, we present a gridless problem accounting for the independence  of fluctuations.
\end{abstract}


\section{Introduction}
\label{sec:intro}

    In super-resolution for biomedical imaging, one wants to retrieve some fine scales details to better study biological structures.
    Indeed, the observed bodies are generally smaller than the Rayleigh limit at \SI{200}{nm}, a length at which the phenomenon of light diffraction comes into play. This diffraction causes a blurring of the image, which can be modelled as a convolution of the image by the PSF. Hence, we want to perform a \emph{deconvolution} \ie\ remove the blur of diffraction to obtain a super-resolved image.
    
    In order to enhance spatial resolution over standard diffraction-limited microscopy techniques and allow imaging of biological structures below the Rayleigh criterion, one can use several techniques of fluorescence microscopy such as STimulated Emission Depletion (STED) \cite{Hell1994} which is harmful implies special equipment and is harmful for the sample, SIM (Structured Illumination Microscopy) \cite{Gustafsson2000} which has a rather limited spatial resolution and Single Molecule Localization Microscopy (SMLM) \cite{Sage2015} which has a low temporal resolution.
    
    Several super-resolution fluorescence methods were proposed to overcome these drawbacks in the last onward decade. An interesting idea lies in the independent stochastic temporal fluctuations of the fluorescent dyes attached to the studied bodies.
    Super-Resolution Optical Fluctuation Imaging (SOFI) \cite{Dertinger2010} exploits the acquired stack by computing its high order $(n)$ statistics that shrinks the point spread function (PSF) by a $\sqrt{n}$ factor. This approach is also less constraining since it can be performed with classic microscope system and conventional dyes. Despite an improvement of temporal resolution, its spatial resolution is a bit restricted compared to SMLM. A similar limitation is noticeable in Super-Resolution Radial Fluctuations (SRRF) \cite{Gustafsson2016} microscopy where local radial symmetries arising from the radial symmetry of the PSF are computed in each frame; or in Sparsity Based Super-Resolution Correlation Microscopy (SPARCOM) \cite{Solomon2019} which exploits both sparsity nature of the fluorophores distribution and statistical prior of uncorrelated blinking.
    
    These techniques are developed within the discrete setting, meaning that the optimisation is performed over a super-resolved image defined on a finer grid than the observed images.
    On the contrary, we define in this paper a method based on variational off-the-grid optimisation, which relies on the definition of a functional which gains convexity at the expense of complexity of the space of optimisation which is the space of Radon measures \cite{Bredies2012, Candes2012}. It allows sparse reconstruction with great precision on both amplitudes and positions of the fluorescent molecules, and sharp theoretical results on the recovery of the spikes under noise regime \cite{Duval2014}.
    
    In this paper, we propose the definition of an operator $\lambdop$ and an off-the-grid variational problem \eqref{eq:covenant} to exploit the independence of the blinking; then we present the super-resolution performed by our method on both simulated and experimental data. It yields compelling results with only one tuning parameter, no learning dataset; not to mention it is suitable for live-cell imaging, common microscopes and conventional dyes.
    To the best knowledge of the authors, this is the first application of the off-the-grid approach to this temporal stochastic fluctuations' context.

\section{Formulation}
\label{sec:theory}

    Temporal quantities, operators and Radon measures will be denoted in bold. The Radon measure $a \deltaop_x$  is modelling a fluorescent molecule and is referred as the spike of amplitude/luminosity $a$ and position $x$.

    \subsection{Off-the-grid inverse problem settings}
    
        $\xx$ is the ambient space where the positions of the spikes live, and is supposed to be a compact of $\R^2$. Spikes are modelled by Dirac measure $\deltaop$: let $\mzer \eqdef \sum_{i=1}^N a_{0,i} \deltaop_{x_{0,i}}$ be the source measure with amplitudes ${a_0} \in \R_+^N$ and positions ${x_0} \in \xx^N$, the problem is to recover this measure from the acquisition $y \in \hilb$ with $\hilb$ a Hilbert space typically $\Ldx$, defined by $y \eqdef \phiop \mzer + w$ where $w \in \hilb$ is typically white Gaussian noise and $\phiop : \radon \to \hilb$ is the PSF defined over the space of Radon measures $\radon$ \cite{Duval2014}:
        
        \begin{align}
            \phiop \m = \int_\xx \varphi(x) \ud \m(x).
        \end{align}
        
        where $\varphi : \xx \to \hilb$ is the measurement kernel, typically the Gaussian approximation 
        \begin{align*}
            \varphi(x) \eqdef   \left(s \mapsto h(s-x) \eqdef \frac{1}{(\sigma \sqrt{2 \pi})^2 } e^{(s-x)^2/ \sigma^2} \right) \in \Ldx
        \end{align*}
        
        where $\sigma$ is the standard deviation of the PSF $h$, followed by a downsampling over $P$ pixels $\xx_i$ namely $\forall i \in \{1, \dots, P\}$: $[\varphi(x)]_i = \int_{\xx_i} h(s-x) \ud s$ .
        To tackle this problem we use the following convex functional \cite{Bredies2012, Azais2015} also called the BLASSO, which stands for Beurling-LASSO, over the infinite dimensional and non-reflexive Banach space $\radon$:
        
        \begin{align}
            \argmin_{\m \in \radon} \dfrac{1}{2} \norm{y - \phiop(\m)}^2_{\mathcal H} + \lambda \mtv
            \tag{$\mathcal P_\lambda (y)$}
            \label{eq:blasso-bruits}
        \end{align}
        
        with regularisation parameter $\lambda > 0$ which accounts for the trade-off between fidelity and sparsity of the reconstruction. The quantity $\mtv$ is the total variation of a measure $\m$, not to be mistaken for the total variation in image denoising. It is the continuous counterpart of the discrete $\ell_1$ norm of vectors:
        \begin{align*}
            \mtv \eqdef \sup \left( {\int_\xx f \ud \m}, f \in \czer, \normsupremum{f} \leq 1 \right).
        \end{align*}

        \subsection{The dynamic off-the-grid extension}
        
        The former section discussed of deconvolution, but in order to reach compelling results and perform super-resolution we need to use the prior on amplitudes fluctuations independence.
        The off-the-grid approach needs some refinements to leverage the temporal context. Let $T > 0$ be the time acquisition.
        Suppose that we have acquisitions, \ie\ images in $\Ldx$, ranging on $[0,T]$: $\y_t \in \Ldx$ is the image frame acquired at time $t \in [0,T]$, and $\y$ is the stack of images with $\y \in \Ldxt$ the space of (Böchner integrable) temporal map \cite{Bredies2019a}. 
        We propose to reconstruct a \emph{dynamic measure} $\mud \in \Ldxm$, temporal (Böchner integrable) map valued on $\radon$, which produces the stack $\y$ by $\y(t) = \phiop \mud(t)$ for $t \in [0,T]$ \aee, in this context its is a map of the form:

        \begin{align}
            t \mapsto \mud(t) \eqdef \sum_{i=1}^N a_i(t) \deltaop_{x_i}.
            \label{eq:mu}
        \end{align}
        
        Note that the positions in the ground-truth dynamic measure $\mud$ do not depend on the time, as we assume that the temporal sample do not move during the acquisition time. One would want an off-the-grid functional similar to \eqref{eq:blasso-bruits} to exploit $\y$ and the statistical independency prior.

    \subsection{Sparse regularisation through cumulants}
        The cumulants can help us find the positions $x_i$ of the ground-truth molecules: let $\ybar \eqdef \int_0^T \y_t \ud t$ be the temporal mean, let $\cumul \in \Ldxd$ be the cumulant of order 2 \ie\ the covariance. Thanks to fluctuations independence \cite{Dertinger2010}, for $u,v \in \xx$ \aee:
    
        \begin{align*}
            \cumul(u,v) &\eqdef  \int_0^T \left(\y_t(u) - \ybar(u) \right) \left(\y_t(v) - \ybar(v) \right) \ud t \\
                &= \textstyle{\sum_{1 \leq i,j}^N \left[ \int_0^T a_i a_j  - \int_0^T a_i \int_0^T a_j \right]} \\
                & \null\qquad \qquad  \times h(u - x_i) h(v - x_j)\\
                &= \textstyle{\sum_{i=1}^N {\left[ \int_0^T a_i^2 - \left( \int_0^T a_i  \right)^2 \right]} h(u - x_i) h(v - x_i)} \\
                &= \sum_{i=1}^N { \underbrace{ M_i}_{\text{Variance of } a_i}  h(u - { x_i}) h(v - { x_i})} \\
                &= \int_\xx h(u - s) h(v - s) \, \dd \mcovar (s).
        \end{align*}
        
        We recognise the action of an operator on a stationary discrete Dirac measure $\mcovar \eqdef \sum_{i=1}^N M_i \deltaop_{x_i}$. From a spatio-temporal acquisition $\y$ we can compute its covariance $\cumul$, which can also be thought as the result of an operator: let us note $\lambdop : \m \longmapsto \cumul$ this operator on covariance.
    
        Since the acquisition kernel of $\phiop$ is $\varphi(x) = h(\cdot -x)$, we can define in the same fashion the kernel of $\lambdop$ by $\phi (x) = h(\cdot -x) h(\cdot -x)$.
        Our operator then writes down for $\m \in \radon$ and $u, v \in \xx$ a.e.: 
        
        \begin{align*}
            \lambdop (\m) (u,v) \eqdef \int_\xx h(u-x) h(v-x) \ud \m (x).
        \end{align*}
        
        Let $\lambda > 0$ be the regularisation parameter, we can introduce our proposed functional to reconstruct from covariance:
        \begin{align}
            \argmin_{\m \in \radon} \dfrac{1}{2} \norm{\cumul - \lambdop(\m)}_{\Ldxd}^2 + \lambda \mtv.
            \tag{$\vb*{\mathcal Q_\lambda (\y)}$}
            \label{eq:covenant}
        \end{align}

        The problem \eqref{eq:covenant} is formally the limit of its discrete counterpart \colorme\ \cite{Stergiopoulou2021} for $\ell_1$-norm on a finer and finer grid. Authors of \cite{Bendory2014} suggest that we can reconstruct spikes close to at most $\sigma / \sqrt{2}$ for \eqref{eq:covenant}, compared to $\sigma$ for the BLASSO on the temporal mean $\ybar$.

\section{Numerical implementation}
\label{sec:numeric}
    
    The BLASSO \ref{eq:blasso-bruits} is an optimisation over the set of Radon measures, an infinite dimensional space which lack algebraic properties. However, there exists multiple algorithms \cite{Candes2012, Denoyelle2018, Chizat2019} to exploit the continuous setting, such as the conditional gradient method also called the Frank-Wolfe (FW) algorithm: if applied to the BLASSO \eqref{eq:blasso-bruits}, it consists in adding iteratively new Dirac masses to the estimated measure. Hence, the interest of this algorithm lies in the fact that it uses only the directional derivatives of the objective functional and that it does not require any Hilbertian structure, contrary to a proximal algorithm formulated in terms of Euclidean distance.

    We use an improved version of the Frank-Wolfe algorithm called the \sfw\ \cite{Denoyelle2018} to further decrease the objective functional, it is presented in Algorithm \ref{algo:sfw} for the case of \eqref{eq:blasso-bruits}. It differs from the original algorithm by its last step, which is a non-convex optimisation in both amplitudes and positions.
    
    \begin{algorithm}[ht!]
        \footnotesize
        \SetAlgoLined
        Initialisation : $m^{[0]} = 0$, $N^{[k]} = 0$.
        
        \For{$k$, $0 \leq k \leq K$}{
            For $m^{[k]} = \sum_{i=1}^{N^{[k]}} a_i^{[k]} \deltaop_{x_i^{[k]}}$ such that $a_i^{[k]} \in \R$, $x_i^{[k]} \in \xx$, find $x_\ast^{[k]} \in \xx$ such that :
                \begin{align*}
                    x_\ast^{[k]} \in \argmax_{x \in \xx} \abs{\eta^{[k]}(x)}, \, \eta^{[k]}(x) \eqdef \dfrac{1}{\lambda} \opadj (\phiop m^{[k]} - y).
                \end{align*}
            \eIf{$\abs{\eta^{[k]} \left(x_\ast^{[k]} \right)} \leq 1$} {
                $m^{[k]}$ is the solution of the BLASSO. Stop.} {
                Let $M^{[k]} \eqdef N^{[k]}+1$, compute $m^{[k + 1/2]} \eqdef \sum_{i=1}^{N^{[k]}} a_i^{[k + 1/2]} \deltaop_{x_i^{[k + 1/2]}} + a^{[k + 1/2]}_{M^{[k]}} \deltaop_{x_\ast^{[k]}}^{[k + 1/2]}$ such that:
                    \begin{align*}
                        a_i^{[k + 1/2]} \in& \argmin_{a \in \R^{M^{[k]}}} \dfrac{1}{2} \norm{y - \phiop_{x^{[k + 1/2]}} (a)}^2_\hilb + \lambda \norm{a}_1 \quad \\ 
                        &\text{for} \quad x^{[k + 1/2]} \eqdef \left(x_1^{[k]},\dots,x_{N^{[k]}}^{[k]}, x_\ast^{[k]} \right)\\
                        & \phiop_{x} (a) \eqdef \sum_{i=1}^N a_i \varphi(x_i), \, x \in \xx^N, a \in \R^N_+.
                    \end{align*}
                    
               Compute $m^{[k+1]} \eqdef \sum_{i=1}^{M^{[k]}} a_i^{[k+1]} \deltaop_{x_i^{[k+1]}}$ such that $\left(a^{[k+1]}, x^{[k+1]} \right) \in$
                \begin{align*}
                    \argmax_{(a,x) \in \R^{M^{[k]}} \times \xx^{M^{[k]}} } \dfrac{1}{2}  \norm{y - \phiop_{x^{\left[k + 1/2\right]}}(a)}^2_\hilb + \lambda \norm{a}_1 .
                \end{align*}
            }
        }
        \caption{\sfw\ for \eqref{eq:blasso-bruits}.}
        \label{algo:sfw}
    \end{algorithm}

    The machinery behind the algorithm can be broke down into 3 operations, namely at the $k$-th step:
    \begin{enumerate}
        \item from the residual, compute $x_\ast^{[k]}$ the position of the new spike estimated at this step;
        \item optimise through LASSO the amplitudes of the reconstructed spikes \ie\ all the spikes reconstructed at the previous steps + the new spike at $x_\ast^{[k]}$;
        \item optimise both amplitudes and positions of all the computed spikes.
    \end{enumerate}
    It can be applied straightforwardly to  \eqref{eq:covenant} by turning the stationary input $y$ into $\cumul$ and the operator $\phiop$ into $\lambdop$.

\section{Results}
\label{sec:results}
    
    \subsection{Simulated data}
        Let us incorporate noise relevant to the physical optic context to the model of acquisition, it is then for $t \in [0,T]$ \aee: $\y_t \eqdef \mathcal{P}(\phiop \mud + \fond) + \w(t)$ where $\mathcal{P}$ is Poisson multiplicative noise (not taken into account in our reconstruction model), $\w$ is white Gaussian noise and $\fond$ is background noise. 
    
\begin{figure}[ht]
        \begin{subfigure}{0.49\columnwidth}
          \centering
            \begin{tikzpicture}[spy using outlines={circle, yellow, magnification=2, size=2.0cm, connect spies}]
                \node{\includegraphics[width=0.95\columnwidth]{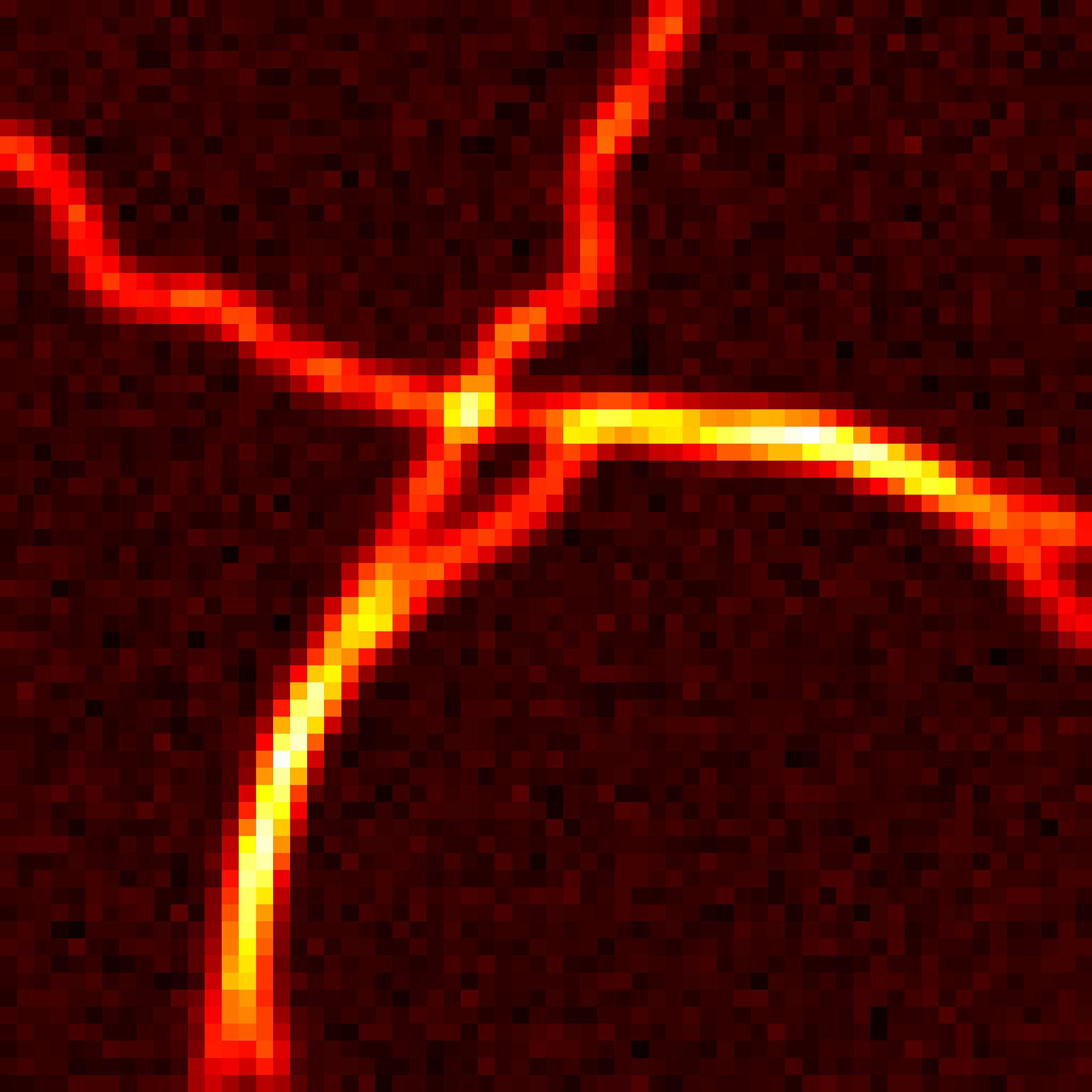}};
                \spy on (0.65,0.4) in node [left] at (1.8, -0.9);
            \end{tikzpicture}
          \caption{}
        \end{subfigure}
        \hfill
        \begin{subfigure}{0.49\columnwidth}
          \centering
            \begin{tikzpicture}[spy using outlines={circle, yellow, magnification=2, size=2.0cm, connect spies}]
                \node{\includegraphics[width=0.95\columnwidth]{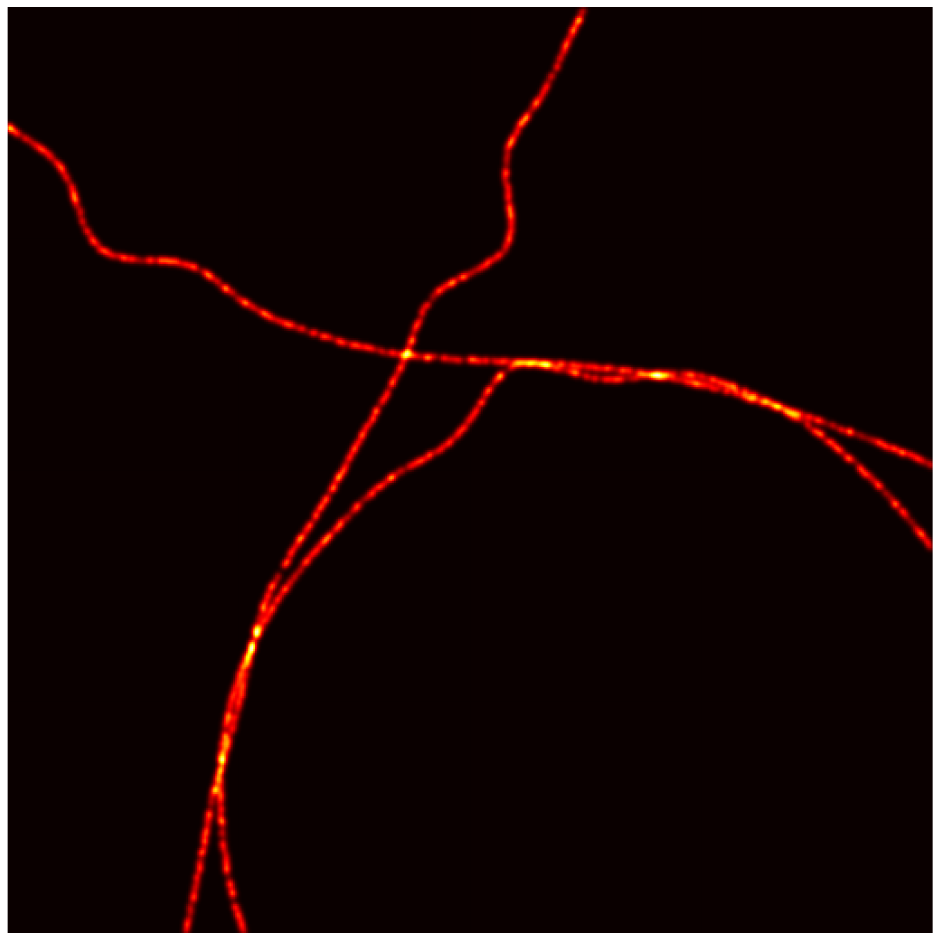}};
                \spy on (0.65,0.4) in node [left] at (1.8, -0.9);
            \end{tikzpicture}
          \caption{}
        \end{subfigure}
        \medskip
        \begin{subfigure}{0.49\columnwidth}
          \centering
            \begin{tikzpicture}[spy using outlines={circle, yellow, magnification=2, size=2.0cm, connect spies}]
                \node {\includegraphics[width=0.95\columnwidth]{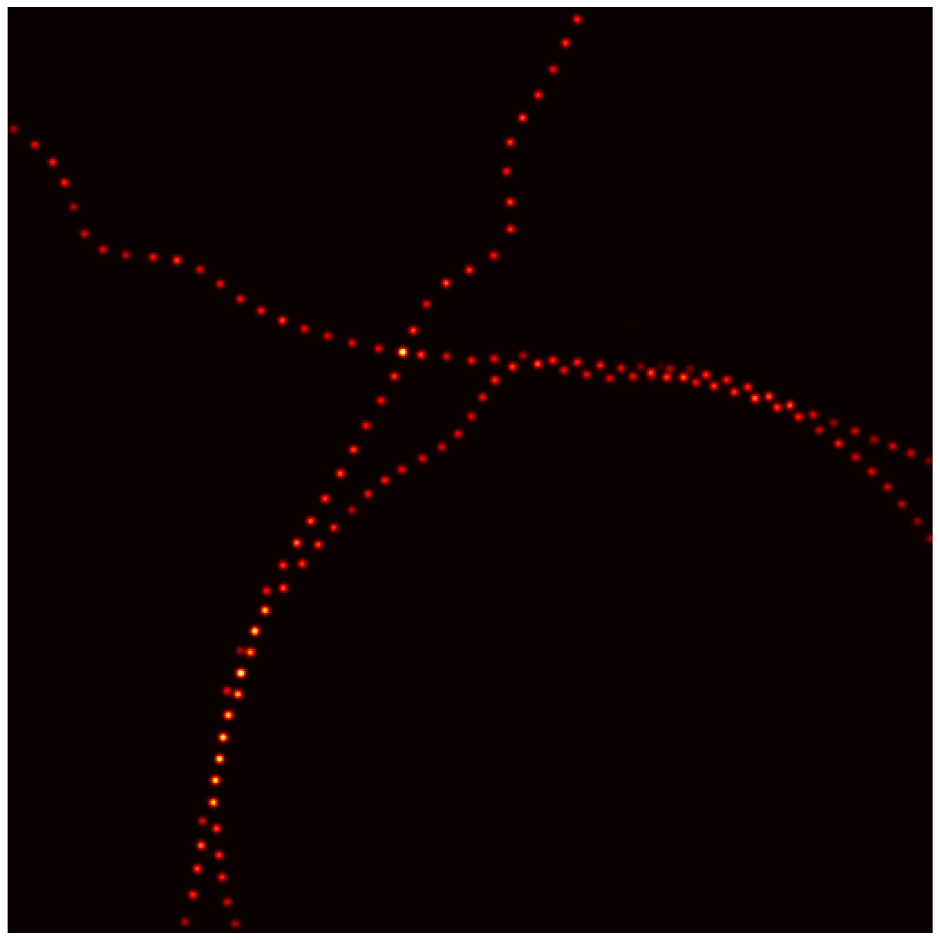}};
                \spy on (0.65,0.4) in node [left] at (1.8, -0.9);
            \end{tikzpicture}
          \caption{}
          \label{fig:sub-third}
        \end{subfigure}
        \begin{subfigure}{0.49\columnwidth}
          \centering
            \begin{tikzpicture}[spy using outlines={circle, yellow, magnification=2, size=2.0cm, connect spies}]
                \node {\includegraphics[width=0.95\columnwidth]{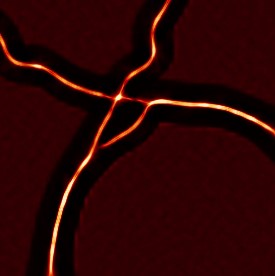}};
                \spy on (0.65,0.4) in node [left] at (1.8, -0.9);
            \end{tikzpicture}
          \caption{}
        \end{subfigure}
        
        \caption{Results for the simulated stack with $T=1000$: (a) first excerpt of the stack $\y$ (b) ground-truth (c) Reconstruction from \eqref{eq:covenant} (d) SRRF.}
        \label{fig:results-filaments}
        \end{figure}

        A first test for our method is made onto images of tubular structure taken from MT0 microtubules training dataset of the SMLM Challenge 2016\footnote{https://srm.epfl.ch/Datasets}, simulating the acquisition process with dyes.
       We generated a stack of $T = 1000$ frames of $64 \times 64$ and 100 frames per second simulated by \texttt{SOFItool} \cite{Girsault2016}, pertaining 8700 molecules covering the tubulins. The FWHM of the PSF is at \SI{229}{nm}, the emitter density is equal to 10.7 emitters/pixel/frame, and the stochastic parameters of fluctuations in \texttt{SOFItool} are chosen at \SI{20}{ms} for on-state average lifetime, \SI{40}{ms} for off-state average lifetime and \SI{20}{s} for average time until bleaching \cite{Stergiopoulou2021}. The stack includes in average 1000 photons/frame per emitting molecule and 100 photons/frame per pixel to simulate the out-of-focus molecules, considered as  background noise. Gaussian noise is also simulated at \SI{20}{db}; to sum-up SNR $\approx$ 10 \si{db}. The results are plotted in Figure \ref{fig:results-filaments}.
    
        One can observe that on the contrary to SRRF, our method is able to recover fine scale details on the interweaving of the tubulins despite the noise level, which is quite hard to achieve by existing methods. Moreover, it does not create artefacts, fits well to the ground-truth inner structure and requires only one tuning parameter \ie\ the $\lambda$.
        Eventually and on the contrary to deep-learning approach, our algorithm does not need any training data set and can be then applied straightforwardly to any stack with known forward operator.

    \subsection{Experimental data}
        
\begin{figure}[ht]
        \begin{subfigure}{0.49\columnwidth}
          \centering
            \begin{tikzpicture}[spy using outlines={circle, yellow, magnification=2, size=2.0cm, connect spies}]
                \node{\includegraphics[width=0.95\columnwidth]{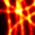}};
            \end{tikzpicture}
          \caption{}
        \end{subfigure}
        \begin{subfigure}{0.49\columnwidth}
          \centering
            \begin{tikzpicture}[spy using outlines={circle, yellow, magnification=2, size=2.0cm, connect spies}]
                \node {\includegraphics[width=0.95\columnwidth]{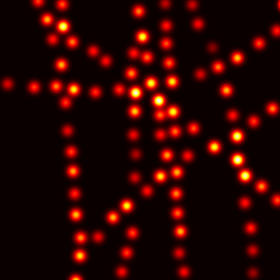}};
            \end{tikzpicture}
          \caption{}
        \end{subfigure}
        \medskip
        \begin{subfigure}{0.49\columnwidth}
          \centering
            \begin{tikzpicture}[spy using outlines={circle, yellow, magnification=2, size=2.0cm, connect spies}]
                \node {\includegraphics[width=0.95\columnwidth]{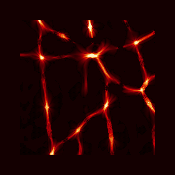}};
            \end{tikzpicture}
          \caption{}
        \end{subfigure}
        
        \caption{Results for real data with $T=500$: (a) mean of the stack $\ybar$ (b) \eqref{eq:covenant} (c) SRRF.}
        \label{fig:results-real-data}
        \end{figure}

        As with SPARCOM \cite{Gustafsson2016} we use high-density stack of SMLM data set: indeed, despite dyes in SMLM do not exhibit blinking but rather an \emph{on} and \emph{off} behaviour, we can make this approximation by noting that there are several molecules per pixel, which results in blinking in each pixel.
        It is thus possible to consider the acquired SMLM stack as a temporal series of blinking.
        The acquired stack is composed of $T=500$ images, blurred by a PSF of FWHM of \SI{351.8}{nm}. The results are plotted in Figure \ref{fig:results-real-data}. 
    
        The reconstruction seems a bit dotted, due to the greater FWHM, but can still recover the structure of the interweaving. tubulins. One can note that we reconstruct spikes close to at most $\sigma/\sqrt{2}$ for \eqref{eq:covenant}.

\section{Conclusion}
\label{sec:conclusion}

We proposed a super-resolution method based on off-the-grid variational approach and blinking fluorescence microscopy, suitable for live-cell imaging and common microscopes. Thanks to the independence of stochastic fluctuations and thanks to the gridless setting, our method can precisely locate the dyes by recasting the problem on the covariance domain. The method succeeds well in the recovery of fine scale details and of hidden structures of the studied interweaving filaments.
Further works plan to better account for the structure of the biological bodies samples in order to reconstruct curves rather than dotted approximation of the filaments curve, and to investigate implementation in 3D \eg\ with bi-plan acquisition. 

\section{Compliance with ethical standards}
\label{sec:ethics}

{
    This work was conducted using biological data available in open access by EPFL SMLM datasets. Ethical approval was not required as confirmed by the license attached with the open access data.
}

\section{Acknowledgments}
\label{sec:acknowledgments}

{
    The work of BL has been supported by the French government, through the UCA DS4H Investments in the Future project managed by the National Research Agency (ANR) with the reference number ANR-17-EURE-0004. The work of LBF has been supported by the French government, through the 3IA Cote d’Azur Investments in the Future project managed by the National Research Agency (ANR) with the reference number ANR-19-P3IA-0002.
}

{

}

\end{document}